\documentclass[11pt]{amsart}
\usepackage{mathrsfs}
\usepackage{amssymb,latexsym}
\usepackage{pstricks,color}
\usepackage{mma}

 \numberwithin{equation}{section}

\newtheorem{theorem}{Theorem}[section]
\newtheorem{proposition}[theorem]{Proposition}
\newtheorem{corollary}[theorem]{Corollary}

\newtheorem{remark}[theorem]{Remark}

\newtheorem{lemma}[theorem]{Lemma}

\def\qed{\hfill $\Box$}
\def\pf{\noindent {\it Proof.} }

\title{Combinatorial identities related to $2\times 2$ submatrices of recursive matrices}

\begin{document}
\maketitle
\begin{center}
Fangfang Cai$^{1}$, Qing-Hu Hou$^{2}$, Yidong Sun\footnote{Corresponding author: Yidong Sun.} and Arthur L.B. Yang$^3$

$^1$School of Science, Dalian Maritime University, 116026 Dalian, P.R. China\\[5pt]
$^2$School of Mathematics, Tianjin University, 300350 Tianjin, P.R. China\\[5pt]
$^3$Center for Combinatorics, Nankai University, 300071 Tianjin, P.R. China\\[5pt]

{\it Emails: cff@dlmu.edu.cn,  qh$\_$hou@tju.edu.cn,  sydmath@dlmu.edu.cn,  yang@nankai.edu.cn }

\end{center}

\subsection*{Abstract}

Recursive matrices are ubiquitous in combinatorics, which have been extensively studied. We focus on the study of the sums of $2\times 2$ minors of certain recursive matrices, the alternating sums of their $2\times 2$ minors, and the sums of their $2\times 2$ permanents. We obtain some combinatorial identities related to these sums, which generalized the work of Sun and Ma in [{\it Electron. J. Combin. 2014}] and [{\it European J. Combin. 2014}].
With the help of the computer algebra package {\tt HolonomicFunctions}, we further get some new identities involving Narayana polynomials.

\medskip

{\bf Keywords}: Recursive matrices; Riordan array; Narayana polynomial; Little
Schr\"{o}der number; Residue.

\noindent {\bf AMS Classification 2010}: 05A10, 05A15, 05A19

{\bf \section{ Introduction } }

Shapiro's Catalan triangle, defined by $\mathcal{B}=(B_{n,k})_{n\geq k\geq 0}$ with
$B_{n,k}=\frac{k+1}{n+1}\binom{2n+2}{n-k}$ \cite[A039598]{Sloane}, appears in various combinatorial settings \cite{Callan, ShapA, ShapB, ShapGet} and has
been paid a lot attention in combinatorics and number theory \cite{ChenChu, Deng, GuoZeng, Gutierrez, Miana, MianaRom, SunFei, SunLu, ZhangPang}.
Table 1.1 illustrates this triangle for $n$ and $k$ up to $5$.
\begin{center}
\begin{eqnarray*}
\begin{array}{|c|cccccc|}\hline
n/k & 0   & 1    & 2    & 3    & 4    & 5       \\\hline
  0 & 1   &      &      &      &      &        \\
  1 & 2   & 1    &      &      &      &         \\
  2 & 5   & 4    & 1    &      &      &         \\
  3 & 14  & 14   & 6    &  1   &      &         \\
  4 & 42  & 48   & 27   &  8   &  1   &        \\
  5 & 132 & 165  & 110  &  44  &  10  &  1     \\\hline
\end{array}
\end{eqnarray*}
Table 1.1. The first values of $B_{n,k}$.
\end{center}

Define another infinite lower triangle $\mathcal{X}=(X_{n,k})_{n\geq k\geq 0}$, based on the triangle $\mathcal{B}$ as follows,
\begin{eqnarray*}
X_{n,k} \hskip-.22cm &=&\hskip-.22cm  \rm{det}\left(\begin{array}{cc}
B_{n, k}   & B_{n,k+1}  \\[5pt]
B_{n+1,k}  & B_{n+1,k+1}
\end{array}\right).
\end{eqnarray*}
Table 1.2 illustrates the triangle $\mathcal{X}$ for $n$ and $k$ up to $4$.
Sun and Ma \cite{SunLu} observed a curious property of the triangle $\mathcal{X}$, namely, the row sums have close relation with the first column of the triangle $\mathcal{B}$.
\begin{center}
\begin{eqnarray*}
\begin{array}{|c|ccccc|c|c|}\hline
n/k & 0   & 1   & 2    & 3    & 4     & row\ sums   \\\hline
  0 & 1   &     &      &      &       &  1=1^2      \\
  1 & 3   & 1   &      &      &       &  4=2^2      \\
  2 & 14  & 10  & 1    &      &       &  25=5^2     \\
  3 & 84  & 90  & 21   & 1    &       &  196=14^2   \\
  4 & 594 & 825 & 308  & 36   &  1    &  1764=42^2  \\\hline
\end{array}
\end{eqnarray*}
Table 1.2. The first values of $X_{n,k}$ and the row sums.
\end{center}
%\begin{center}
%\begin{eqnarray*}
%\begin{array}{|c|ccccc|c|c|}\hline
%n/k & 0   & 1   & 2    & 3    & 4     & row\ sums   & alternating\ row\ sums   \\\hline
%  0 & 1   &     &      &      &       &  1=1^2      &  1                      \\
%  1 & 3   & 1   &      &      &       &  4=2^2      &  2                      \\
%  2 & 14  & 10  & 1    &      &       &  25=5^2     &  5                       \\
%  3 & 84  & 90  & 21   & 1    &       &  196=14^2   &  14                     \\
%  4 & 594 & 825 & 308  & 36   &  1    &  1764=42^2  &  42                      \\\hline
%\end{array}
%\end{eqnarray*}
%Table 1.2. The first values of $X_{n,k}$ and the (alternating) row sums.
%\end{center}

Note that Shapiro's Catalan triangle is a special class of recursive matrices, which were introduced by Aigner \cite{Aigner} in the study of Catalan-like numbers. Hence it is natural to consider whether similar phenomena would happen to other recursive matrices.
Following Aigner, we say that $\mathcal{A}$ is a recursive matrix if $\mathcal{A}=\mathcal{A}^{\sigma,\tau}=(A_{n,k})_{n\geq k\geq 0}$ for some pair of sequences $\sigma=(\sigma_0,\sigma_1,\sigma_2,\ldots)$ and $\tau=(\tau_1,\tau_2,\tau_3,\ldots)$ with $\tau_i\neq 0$ for each $i\geq 1$, where
\begin{align}\label{eq-rec-matrix}
A_{0,0}&=1,\quad A_{0,k}=0, \ (k>0),\\
A_{n,k}&=A_{n-1,k-1}+\sigma_k A_{n-1,k}+\tau_{k+1} A_{n-1,k+1}, \ (n\geq 1).\nonumber
\end{align}
It is well known that each entry $A_{n,k}$ admits a weighted partial Motzkin path interpretation, see \cite{Aigner, Flajolet, Pan}. Using the weighted partial Motzkin paths,
Sun and Ma \cite{SunLu} obtained the following general result.

\begin{theorem}\label{thm-1-1}
For any integers $n, r\geq 0$ and $m\geq \ell\geq 0$, set $M_r=\min\{n+r+1, m+r-\ell\}$. If $(A_{n,k})_{n\geq k\geq 0}$ is the recursive matrix $\mathcal{A}^{\sigma,\tau}$ corresponding to $\sigma=(x,y,y,\ldots)$ and $\tau=(1,1,1,\ldots)$, namely,
$\sigma_0=x,\sigma_i=y,\tau_i=1$ for any $i\geq 1$, then we have
\begin{align}
\sum_{k=0}^{M_r}
\det\left(\begin{array}{cc}
A_{n,k}     & A_{m,k+\ell+1} \\[5pt]
A_{n+r+1,k} & A_{m+r+1,k+\ell+1}
\end{array}\right)
&=\sum_{i=0}^{r} A_{n+i,0}A^{*}_{m+r-i,\ell}. \label{eqn 1.1}
\end{align}
where $A^{*}_{n,k}$ denotes $A_{n,k}$ in the $x=y$ case.
\end{theorem}

For $x=y=2$, $m=n$ and $r=\ell=0$, the above theorem reduces to the aforementioned observation for Shapiro's Catalan triangle.
Sun and Ma \cite{SunFei} also obtained a similar result for the permanent case.

\begin{theorem}\label{thm-1-2}
If $(A_{n,k})_{n\geq k\geq 0}$ is the recursive matrix
$\mathcal{A}^{\sigma,\tau}$ corresponding to $\sigma=(y,y,y,\ldots)$
and $\tau=(1,1,1,\ldots)$, then for any integers $m, n, r$ with
$m\geq n\geq 0$ we have
\begin{align} \label{eqn 1.3}
\sum_{k=0}^{m}\rm{per}\left(\begin{array}{cc}
A_{n,k} & A_{n+r,k+1} \\[5pt]
A_{m,k} & A_{m+r,k+1}
\end{array}\right)
&= A_{m+n+r,1}+H_{n,m}(r),
\end{align}
where ${\rm per}(A)$ denotes the permanent of a square matrix $A$, and
\begin{eqnarray*}
H_{n,m}(r)=\left\{\begin{array}{ll}
\sum_{i=0}^{r-1}A_{n+i,0}A_{m+r-i-1,0},    & {\rm{if}} \ r\geq 1, \\[5pt]
0,                                         & {\rm{if}} \ r=0,      \\[5pt]
-\sum_{i=1}^{|r|}A_{n-i,0}A_{m-|r|+i-1,0}, & {\rm{if}} \ r\leq
-1.
\end{array}\right.
\end{eqnarray*}
\end{theorem}

When the parameters $(x,y)$ are specialized, Theorems \ref{thm-1-1} and \ref{thm-1-2} can produce many combinatorial identities related to
Catalan numbers, Motzkin numbers and other combinatorial sequences \cite{SunFei, SunLu}. Despite this,
there still have many classical combinatorial triangles do not fall into the above framework. For example, if we consider the triangular
array $\mathcal{S}=(s_{n,k})_{n\geq k\geq 0}$ \cite[A110440]{Sloane} formed by the little
Schr\"{o}der numbers by the following recurrence
\begin{eqnarray} \label{eqn 1.4}
s_{n,k}=s_{n-1,k-1}+3s_{n-1,k}+2s_{n-1,k+1},
\end{eqnarray}
with $s_{n,n}=1$ for $n\geq 0$ and $s_{n,k}=0$ for $n<k$ or $n,k<0$.
Table 1.3 illustrates the triangle $\mathcal{S}$ for $n$ and
$k$ up to $5$. The leftmost column $s(n,0)$ is the sequence of the
little Schr\"{o}der numbers $s_{n}$ \cite[A001003]{Sloane}.

\begin{center}
\begin{eqnarray*}
\begin{array}{|c|cccccc|}\hline
n/k & 0   & 1    & 2    & 3    & 4    & 5       \\\hline
  0 & 1   &      &      &      &      &        \\
  1 & 3   & 1    &      &      &      &         \\
  2 & 11  & 6    & 1    &      &      &         \\
  3 & 45  & 31   & 9    &  1   &      &         \\
  4 & 197 & 156  & 60   &  12  &  1   &        \\
  5 & 903 & 785  & 360  &  98  &  15  &  1     \\\hline
\end{array}
\end{eqnarray*}
Table 1.3. The first values of $s_{n,k}$.
\end{center}

Let $\mathcal{L}=(L_{n,k})_{n\geq k\geq 0}$ be the infinite lower
triangles defined on the triangle $\mathcal{S}$ by
\begin{eqnarray*}
L_{n,k} \hskip-.22cm &=&\hskip-.22cm \rm{det}\left(\begin{array}{cc}
s_{n, k}   & s_{n,k+1}  \\[5pt]
s_{n+1,k}  & s_{n+1,k+1}
\end{array}\right).
\end{eqnarray*}
Table 1.4 illustrates the triangle $\mathcal{L}$ for small $n$ and
$k$ up to $4$, together with the weighted row sums. It indicates
that the weighted row sums have close relation with the first column
of the triangle $\mathcal{S}$.

\begin{center}
\small\begin{eqnarray*}
\begin{array}{|c|ccccc|r|r|}\hline
n/k & 0   & 1   & 2    & 3    & 4     &  weighted\ row \  sums            \\\hline
  0 & 1   &     &      &      &       &  1\times 1=1                              \\
  1 & 7   & 1   &      &      &       &  7\times 1+1\times 2=3^2                          \\
  2 & 71  & 23  & 1    &      &       &  71\times 1+23\times 2+1\times 4=11^2                      \\
  3 & 913 & 456 & 48   & 1    &       &  913\times 1+456\times 2+48\times 4+1\times 8=45^2                 \\
  4 & 13777 & 9060 & 1560 & 82   &  1    & 13777\times 1+9060\times 2+1560\times 4+82\times 8+1\times 16=197^2     \\\hline
\end{array}
\end{eqnarray*}
Table 1.4. The first values of $L_{n,k}$ and the weighted row sums.
\end{center}

Motivated by the above phenomenon, we further study the weighted sums of $2\times 2$ minors and $2\times 2$
permanents of recursive matrices. We obtain the following results, which generalize Theorem \ref{thm-1-1} and Theorem \ref{thm-1-2} respectively.

\begin{theorem}\label{main-thm-1}
For any integers $n, r\geq 0$ and $m\geq \ell\geq 0$, set $M_r=\min\{n+r+1, m+r-\ell\}$. If $(A_{n,k})_{n\geq k\geq 0}$ is the recursive matrix $\mathcal{A}^{\sigma,\tau}$ corresponding to $\sigma=(x,y,y,\ldots)$ and $\tau=(z,z,z,\ldots)$, then we have
\begin{align}
\sum_{k=0}^{M_r}
z^k\det\left(\begin{array}{cc}
A_{n,k}     & A_{m,k+\ell+1} \\[5pt]
A_{n+r+1,k} & A_{m+r+1,k+\ell+1}
\end{array}\right)
&=\sum_{i=0}^{r} A_{n+i,0}A^{*}_{m+r-i,\ell}. \label{main-eqn-1}
\end{align}
where $A^{*}_{n,k}$ denotes $A_{n,k}$ in the $x=y$ case.
\end{theorem}

\begin{theorem}\label{main-thm-2}
If $(A_{n,k})_{n\geq k\geq 0}$ is the recursive matrix $\mathcal{A}^{\sigma,\tau}$ corresponding to $\sigma=(y,y,y,\ldots)$ and $\tau=(z,z,z,\ldots)$, then for any integers $m, n, r$ with $m\geq n\geq 0$ we have
\begin{align} \label{main-eqn-2}
\sum_{k=0}^{m}z^k\rm{per}\left(\begin{array}{cc}
A_{n,k} & A_{n+r,k+1} \\[5pt]
A_{m,k} & A_{m+r,k+1}
\end{array}\right)
&= A_{m+n+r,1}+H_{n,m}(r),
\end{align}
where
\begin{eqnarray*}
H_{n,m}(r)=\left\{\begin{array}{ll}
\sum_{i=0}^{r-1}A_{n+i,0}A_{m+r-i-1,0},    & {\rm if} \ r\geq 1, \\[5pt]
0,                                         & {\rm if} \ r=0,      \\[5pt]
-\sum_{i=1}^{|r|}A_{n-i,0}A_{m-|r|+i-1,0}, & {\rm if} \ r\leq
-1.
\end{array}\right.
\end{eqnarray*}
\end{theorem}

The remainder of this paper is organized as follows. Section \ref{sect-2} will be devoted to proving Theorems \ref{main-thm-1} and \ref{main-thm-2}. In Section \ref{sect-3}, we study the recursive matrix $\mathcal{A}^{\sigma,\tau}$ corresponding to $\sigma=(z+1,z+1,z+1,\ldots)$ and $\tau=(z,z,z,\ldots)$, and give an expression of some weighted sums of its $2\times 2$ minors in terms of Narayana polynomials. In Section \ref{sect-4}, we evaluate the weighted sums of $2\times 2$ minors of the recursive matrix $\mathcal{A}^{\sigma,\tau}$ corresponding to $\sigma=(x,y,y,\ldots)$ and $\tau=(0,0,0,\ldots)$.

\section{Proofs of Theorems \ref{main-thm-1} and \ref{main-thm-2}} \label{sect-2}

The main objective of this section is to prove Theorems \ref{main-thm-1} and \ref{main-thm-2}. We also give some of their corollaries. Before proving Theorems \ref{main-thm-1} and \ref{main-thm-2}, let us first note the following useful lemma.

\begin{lemma}\label{lemma-key} If $z\neq 0$, then we have
\begin{align}
{A}_{n,k}&=z^{\frac{n-k}{2}}\bar{A}_{n,k},
\end{align}
where $(A_{n,k})_{n\geq k\geq 0}$ is the recursive matrix $\mathcal{A}^{\sigma,\tau}$ corresponding to $\sigma=(x,y,y,\ldots)$ and $\tau=(z,z,z,\ldots)$, and $(\bar{A}_{n,k})_{n\geq k\geq 0}$ is the recursive matrix $\mathcal{\bar{A}}^{\sigma,\tau}$ corresponding to $\sigma=(\frac{x}{\sqrt{z}},\frac{y}{\sqrt{z}},\frac{y}{\sqrt{z}},\ldots)$ and $\tau=(1,1,1,\ldots)$.
\end{lemma}

\pf
By definition, it suffices to show that the matrix $(z^{\frac{n-k}{2}}\bar{A}_{n,k})_{n\geq k\geq 0}$ also satisfies the following recurrence relation:
\begin{align*}
z^{\frac{0}{2}}\bar{A}_{0,0}&=1,\quad z^{\frac{-k}{2}}\bar{A}_{0,k}=0 \ (k>0),\\
\left(z^{\frac{n}{2}}\bar{A}_{n,0}\right)&=x\left(z^{\frac{n-1}{2}}
\bar{A}_{n-1,0}\right)+z\left(z^{\frac{n-2}{2}}\bar{A}_{n-1,1}\right),\ (n\geq 1),\\
\left(z^{\frac{n-k}{2}}\bar{A}_{n,k}\right)&=\left(z^{\frac{n-k}{2}}
\bar{A}_{n-1,k-1}\right)+y\left(z^{\frac{n-k-1}{2}}\bar{A}_{n-1,k}\right)+z
\left(z^{\frac{n-k-2}{2}}\bar{A}_{n-1,k+1}\right), \ (n,k\geq 1).
\end{align*}
Since $z\neq 0$, the above recurrence relation is equivalent to
\begin{align*}
{A}_{0,0}&=1,\quad {A}_{0,k}=0, \ (k>0),\\
{A}_{n,0}&=\frac{x}{\sqrt{z}}{A}_{n-1,0}+{A}_{n-1,1},\ (n\geq 1),\\
{A}_{n,k}&={A}_{n-1,k-1}+\frac{y}{\sqrt{z}}
{A}_{n-1,k}+{A}_{n-1,k+1},\ (n,k\geq 1).
\end{align*}
Again by definition, this is just the recurrence satisfied by the matrix $({A}_{n,k})_{n\geq k\geq 0}$.
\qed

Now we are able to prove Theorem \ref{main-thm-1}.

\allowdisplaybreaks \noindent{Proof of Theorem \ref{main-thm-1}.}
Since both sides of \eqref{main-eqn-1} are polynomials in $z$, it
suffices to prove \eqref{main-eqn-1} for any $z\neq 0$. Let
$A^{*}_{n,k}=A_{n,k}$ and $\bar{A}^{*}_{n,k}=\bar{A}_{n,k}$ when
$x=y$. For the left hand side of \eqref{main-eqn-1}, by Lemma
\ref{lemma-key} we have
\begin{align*}
\mbox{(LHS) of \eqref{main-eqn-1}}
&=\sum_{k=0}^{M_r}
z^k\det\left(\begin{array}{cc}
z^{\frac{n-k}{2}}\bar{A}_{n,k}    & z^{\frac{m-k-\ell-1}{2}}\bar{A}_{m,k+\ell+1} \\[5pt]
z^{\frac{n+r+1-k}{2}}\bar{A}_{n+r+1,k} & z^{\frac{m+r-k-\ell}{2}}\bar{A}_{m+r+1,k+\ell+1}
\end{array}\right)\\[8pt]
&=z^{\frac{n+m+r-\ell}{2}}\sum_{k=0}^{M_r}
\det\left(\begin{array}{cc}
\bar{A}_{n,k}    & \bar{A}_{m,k+\ell+1} \\[5pt]
\bar{A}_{n+r+1,k} & \bar{A}_{m+r+1,k+\ell+1}
\end{array}\right)\\[8pt]
&=z^{\frac{n+m+r-\ell}{2}}  \sum_{i=0}^{r} \bar{A}_{n+i,0}\bar{A}^{*}_{m+r-i,\ell}\hskip0.5cm \mbox{ (by Theorem \ref{thm-1-1})}\\[8pt]
&=\sum_{i=0}^{r} \left(z^{\frac{n+i}{2}}\bar{A}_{n+i,0}\right)\left(z^{\frac{m+r-i-\ell}{2}}\bar{A}^{*}_{m+r-i,\ell}\right)\\[8pt]
&=\sum_{i=0}^{r} {A}_{n+i,0}{A}^{*}_{m+r-i,\ell}\hskip0.5cm \mbox{
(by Lemma \ref{lemma-key})},
\end{align*}
which is just the right hand side of \eqref{main-eqn-1}.
This completes the proof.
\qed

Taking $x=y=3$, $z=2$ and $r=\ell=0$ in Theorem \ref{main-thm-1}, we immediately obtain the following result, which answers the aforementioned phenomenon when $m=n$.

\begin{corollary} For any integers $n,m\geq 0$, there holds
\begin{align*}
\sum_{k=0}^{m}2^k\det\left(
\begin{array}{cc}
s_{n,k}      & s_{m,k+1}       \\[5pt]
s_{n+1,k}    & s_{m+1,k+1}
\end{array}\right)
&= s_{n}s_{m},
\end{align*}
where $s_n$ is the $n$-th little Schr\"{o}der number, and $s_{n,k}$ is defined by \eqref{eqn 1.4}.
\end{corollary}

We proceed to prove Theorem \ref{main-thm-2} in the same way as Theorem \ref{main-thm-1} is proved.

\noindent{Proof of Theorem \ref{main-thm-2}.}
Since both sides of \eqref{main-eqn-2} are polynomials in $z$, it suffices to prove \eqref{main-eqn-2} for any $z\neq 0$.
Let $(\bar{A}_{n,k})_{n\geq k\geq 0}$ denote the recursive matrix $\mathcal{\bar{A}}^{\sigma,\tau}$ corresponding to $\sigma=(\frac{y}{\sqrt{z}},\frac{y}{\sqrt{z}},\frac{y}{\sqrt{z}},\ldots)$ and $\tau=(1,1,1,\ldots)$. Thus by Lemma \ref{lemma-key} we have
\begin{align*}
\mbox{(LHS) of \eqref{main-eqn-2}}
&=\sum_{k=0}^{M_r}
z^k\rm{per}\left(\begin{array}{cc}
z^{\frac{n-k}{2}}\bar{A}_{n,k}    & z^{\frac{n+r-k-1}{2}}\bar{A}_{n+r,k+1}\\[5pt]
z^{\frac{m-k}{2}}\bar{A}_{m,k} & z^{\frac{m+r-k-1}{2}}\bar{A}_{m+r,k+1}
\end{array}\right)\\[8pt]
&=z^{\frac{n+m+r-1}{2}}\sum_{k=0}^{M_r}
\rm{per}\left(\begin{array}{cc}
\bar{A}_{n,k}    & \bar{A}_{n+r,k+1}\\[5pt]
\bar{A}_{m,k} & \bar{A}_{m+r,k+1}
\end{array}\right)\\[8pt]
&=z^{\frac{n+m+r-1}{2}}\left(\bar{A}_{m+n+r,1}+\bar{H}_{n,m}(r)\right) \hskip0.5cm\mbox{ (by Theorem \ref{thm-1-2})},
\end{align*}
where
\begin{eqnarray*}
\bar{H}_{n,m}(r)=\left\{\begin{array}{ll}
\sum_{i=0}^{r-1}\bar{A}_{n+i,0}\bar{A}_{m+r-i-1,0},    & {\rm if} \ r\geq 1, \\[5pt]
0,                                                     & {\rm if} \ r=0,      \\[5pt]
-\sum_{i=1}^{|r|}\bar{A}_{n-i,0}\bar{A}_{m-|r|+i-1,0}, & {\rm if} \ r\leq
-1.
\end{array}\right.
\end{eqnarray*}
By Lemma \ref{lemma-key} it is routine to verify that
\begin{align*}
z^{\frac{n+m+r-1}{2}}\left(\bar{A}_{m+n+r,1}+\bar{H}_{n,m}(r)\right)= {A}_{m+n+r,1}+{H}_{n,m}(r).
\end{align*}
This completes the proof.
\qed

Taking $y=3$, $z=2$ and $r=0$ in Theorem \ref{main-thm-2}, we immediately obtain the following result.
\begin{corollary}
For any integers $m\geq n\geq 0$, there holds
\begin{eqnarray*}
\sum_{k=0}^{m}2^k\rm{per}\left(
\begin{array}{cc}
s_{n,k}      & s_{m,k+1}       \\[5pt]
s_{n+1,k}    & s_{m+1,k+1}
\end{array}\right)
\hskip-.22cm &=&\hskip-.22cm  s_{m+n,1},
\end{eqnarray*}
where $s_{n,k}$ is defined by \eqref{eqn 1.4}.
\end{corollary}

\section{The recursive matrix $\mathcal{A}^{(z+1,z+1,z+1,\ldots),(z,z,z,\ldots)}$} \label{sect-3}

In this section we aim to study certain weighted sums of
$2\times 2$  minors of the recursive matrix $\mathcal{A}^{(z+1,z+1,z+1,\ldots),(z,z,z,\ldots)}$.

Let us first determine the entries of $\mathcal{A}^{(z+1,z+1,z+1,\ldots),(z,z,z,\ldots)}$ explicitly.
To this end, we note a connection between recursive matrices and Riordan arrays, see also \cite{wangyi1,wangyi2}. Recall that a
{\it Riordan array}, denoted $(g(v), f(v))$, is an infinite
lower triangular matrix $\mathcal{R}=(R_{n,k})_{n\geq k\geq 0}$ with nonzero entries on the main diagonal,
such that $R_{n, k}=[v^n]g(v)(f(v))^{k}$ for $n\geq k$, where $g(v)=1+g_1v +g_2v^2+\cdots$ and
$f(v)=f_1v + f_2v^2 + \cdots$ with $f_1\neq 0$ are two formal power series. The Riordan array $\mathcal{R}=(R_{n,k})_{n\geq k\geq 0}=(g(v), f(v))$ can also be characterized by the following recurrence relations
\begin{align*}
R_{0,0}&=1,\quad R_{0,k}=0, \ (k>0),\\
R_{n,0}&=z_0R_{n-1,0}+z_1R_{n-1,1}+z_2R_{n-1,2}+\cdots , \ (n\geq 1),\\
R_{n,k}&=a_0R_{n-1,k-1}+a_1R_{n-1,k}+a_2R_{n-1,k+1}+\cdots , \ (n,k\geq 1),
\end{align*}
and moreover
\begin{align}
g(v)&=\frac{1}{1-vZ(f(v))},\label{eq-gf-1}\\
f(v)&=vA(f(v)),\label{eq-gf-2}
\end{align}
where
\begin{align*}
A(v)&=a_0+a_1v+a_2v^2+\cdots,\\
Z(v)&=z_0+z_1v+z_2v^2+\cdots.
\end{align*}
Thus $\mathcal{R}$ is uniquely determined by $A(v)$ and $Z(v)$.
For more information on Riordan arrays, see \cite{cheon, ShapB, ShapGet, Sprug}. We have the following result.

\begin{lemma} \label{lem-connection} The recursive matrix
$\mathcal{A}^{(x,y,y,\ldots),(z,z,z,\ldots)}$ is a Riordan array
$(g(v),f(v))$, where
\begin{align}
f(v)&=\frac{1-yv-\sqrt{(1-yv)^2-4zv^2}}{2zv}, \label{eqn 2.3} \\[8pt]
% why we do not take the positive part? Since f(0)=0.
g(v)&=\frac{1-2xv+yv-\sqrt{(1-yv)^2-4zv^2}}{2(y-x)(1-xv)v+2zv^2}. \label{eqn 2.4}
% Here $2(y-x)(1-xv) v$ should be correct instead of $2(y-x)(1-xv) z$. Yes.
\end{align}
\end{lemma}

\pf By definition the recursive matrix
$\mathcal{A}^{(x,y,y,\ldots),(z,z,z,\ldots)}=(A_{n,k})_{n\geq k\geq 0}$ satisfies the following recurrence:
\begin{align*}
A_{0,0}&=1,\quad A_{0,k}=0, \ (k>0),\\
A_{n,0}&=xA_{n-1,0}+zA_{n-1,1},\ (n\geq 1),\\
A_{n,k}&=A_{n-1,k-1}+yA_{n-1,k}+zA_{n-1,k+1},\ (n,k\geq 1).
\end{align*}
Comparing the above recurrence with the recurrence relation satisfied by a Riordan array, it is clear that $\mathcal{A}^{(x,y,y,\ldots),(z,z,z,\ldots)}$ is a Riordan array
with
\begin{align*}
A(v)&=1+yv+zv^2,\\
Z(v)&=x+zv.
\end{align*}
By \eqref{eq-gf-2}, we obtain that
\begin{align*}
f(v)&=v(1+yf(v)+zf(v)^2).
\end{align*}
Solving this equation, we get that
\begin{align*}
f(v)&=\frac{1-yv-\sqrt{(1-yv)^2-4zv^2}}{2zv}.
\end{align*}
Substitute $f(v)$ into the right-hand side of \eqref{eq-gf-1} to complete the proof of \eqref{eqn 2.4}. \qed

Note that the above lemma could also be deduced by using the partial weighted Motzkin paths, as was done in \cite{SunLu}. Now we can determine the entries of $\mathcal{A}^{(z+1,z+1,z+1,\ldots),(z,z,z,\ldots)}$.
We have the following result.

\begin{proposition} Let $N_{n,k}(z)$ denote the entry of $\mathcal{A}^{(z+1,z+1,z+1,\ldots),(z,z,z,\ldots)}$ in the $n$-th row and $k$-th column. Then
\begin{align}\label{eq-narayana}
N_{n,k}(z)&=\sum_{i=0}^{n-k}\frac{k+1}{n+1}\binom{n+1}{i}\binom{n+1}{i+k+1}z^{i}. \end{align}
\end{proposition}

\pf By Lemma \ref{lem-connection}, the recursive matrix $\mathcal{A}^{(z+1,z+1,z+1,\ldots),(z,z,z,\ldots)}$ is the Riordan array $(g(v), f(v))$ with
\begin{align}\label{eq-ff}
f(v)=vg(v)=\frac{1-v-zv-\sqrt{1-2v+v^2-2vz-2zv^2+v^2z^2}}{2zv}.
\end{align}
From \eqref{eq-gf-2} it follows that $f(v)$ satisfies
\begin{eqnarray}\label{eqn 2.5}
v(1+f(v))(1+zf(v))=f(v).
\end{eqnarray}
By applying the Lagrange inversion formula, we obtain that
\begin{align*}
 N_{n,k}(z)&=[v^n]g(v)(f(v))^{k}\\
&=[v^{n+1}](f(v))^{k+1}\\
&=[v^{n}]\frac{k+1}{n+1}v^{k}(1+v)^{n+1}(1+zv)^{n+1}\\
 &=\sum_{i=0}^{n-k}\frac{k+1}{n+1}\binom{n+1}{i}\binom{n+1}{i+k+1}z^{i},
\end{align*}
as was shown in Appendix D of \cite{Deut}. This completes the proof.
\qed

Clearly, $N_{n,k}(z)$ is a symmetric polynomial in $z$, i.e. $N_{n,k}(z)=z^{n-k}A_{n,k}(z^{-1})$, and
$N_{n,k}(z)$ in the $k=0$ case leads to the classical Narayana polynomial $N_n(z)$, namely,
\begin{eqnarray*}
N_{n}(z) \hskip-.25cm &=& \hskip-.25cm   \sum_{i=0}^{n}\frac{1}{n+1}\binom{n+1}{i}\binom{n+1}{i+1}z^{i}.
\end{eqnarray*}
The few values of $N_{n,k}(z)$ are illustrated in Table
2.1.
\begin{center}
\small\begin{eqnarray*}
\begin{array}{|c|l|l|l|l|l|l}\hline
n/k & 0 & 1 & 2 & 3 & 4 \\\hline
0 & 1                           &                          &                     &                &        \\
1 & z+1                         & 1                        &                     &                &         \\
2 & z^2+3z+1                    & 2z+2                     & 1                   &                &         \\
3 & z^3+6z^2+6z+1               & 3z^2+8z+3                & 3z+3                & 1              &         \\
4 & z^4+10z^3+20z^2+10z+1       & 4z^3+20z^2+20z+4         & 6z^2+15z+6          & 4z+4           & 1        \\\hline
%5 & z^5+15z^4+50z^3+50z^2+15z+1 & 5z^4+40z^3+75z^2+40z+5   & 10z^3+45z^2+45z+10  & 10z^2+24z+10   & 5z+5  & 1     \\\hline
\end{array}
\end{eqnarray*}
Table 2.1. The first values of $N_{n,k}(z)$.
\end{center}

Note that the triangle $\{N_{n,k}(z)\}_{n\geq k\geq 0}$ in the cases $z=0$, $z=1$ and $z=2$ leads respectively to the Pascal triangle, Shapiro's Catalan
triangle and the triangle $\mathcal{S}$ aforementioned, that is, $N_{n,k}(0)=\binom{n}{k}$, $N_{n,k}(1)=B_{n,k}$ and $N_{n,k}(2)=s_{n,k}$ for $n\geq k\geq 0$.
By using certain weighted NSEW-paths, Cigler \cite{Cigler} gave another combinatorial interpretation, together with another expression for $N_{n,k}(z)$,
\begin{eqnarray*}
N_{n,k}(z) \hskip-.25cm &=& \hskip-.25cm \sum_{i=0}^{[\frac{n-k}{2}]}\frac{k+1}{i+k+1}\binom{n}{2i+k}\binom{k+2i}{i}z^{i}(1+z)^{n-k-2i}.
\end{eqnarray*}

Now the following results immediately follow from Theorems \ref{main-thm-1} and \ref{main-thm-2}.

\begin{corollary} For any integers $n, r\geq 0$ and $m\geq \ell\geq 0$, set $M_r=\min\{n+r+1, m+r-\ell\}$.
Then there holds
\begin{eqnarray*}
\sum_{k=0}^{M_r}z^k\det\left(\begin{array}{cc}
N_{n,k}(z)     & N_{m,k+\ell+1}(z)      \\[5pt]
N_{n+r+1,k}(z) & N_{m+r+1,k+\ell+1}(z)
\end{array}\right)
\hskip-.22cm &=&\hskip-.22cm \sum_{i=0}^{r} N_{n+i,0}(z)N_{m+r-i,\ell}(z).
\end{eqnarray*}
Specially,
\begin{eqnarray*}
\sum_{k=0}^{m}z^k\det\left(\begin{array}{cc}
N_{n,k}(z)     & N_{m,k+1}(z)      \\[5pt]
N_{n+1,k}(z)   & N_{m+1,k+1}(z)
\end{array}\right)
\hskip-.22cm &=&\hskip-.22cm  N_{n}(z)N_{m}(z).
\end{eqnarray*}
\end{corollary}

\begin{corollary}
For any integers $m\geq n\geq 0$, there holds
\begin{eqnarray*}
\sum_{k=0}^{m}z^k\rm{per}\left(\begin{array}{cc}
N_{n,k}(z) & N_{n,k+1}(z) \\[5pt]
N_{m,k}(z) & N_{m,k+1}(z)
\end{array}\right)
\hskip-.22cm &=&\hskip-.22cm N_{m+n,1}(z).
\end{eqnarray*}
\end{corollary}

Next we consider an alternating summation of weighted $2\times 2$ minors of the recursive matrix $\mathcal{A}^{(z+1,z+1,z+1,\ldots),(z,z,z,\ldots)}=(N_{n,k}(z))_{n,k\geq 0}$. We have the following main result.

\begin{theorem}\label{theo 3.6}
For any integers $m, n\geq 0$, let
\begin{eqnarray}\label{eqn 3.5}
F_{m,n}(z) \hskip-.22cm &=&\hskip-.22cm
\sum_{k=0}^{n}(-z)^k\det\left(\begin{array}{cc}
N_{n,k}(z)     & N_{n,k+1}(z)      \\[5pt]
N_{m,k}(z)     & N_{m,k+1}(z)
\end{array}\right).
\end{eqnarray}
Then there holds the recurrence
\begin{eqnarray}\label{eqn 3.6}
F_{m,n}(z)=2(z+1)F_{m-1,n}(z)-F_{m-1,n+1}(z)
\end{eqnarray}
with $F_{n,n}(z)=0$ and $F_{n+1,n}(z)=N_n(z^2)$. Moreover, for $m> n\geq 0$, we have
\begin{eqnarray}\label{eqn 3.7}
F_{m,n}(z)=\sum_{j=0}^{[\frac{m-n-1}{2}]}(-1)^j\binom{m-n-1-j}{j}N_{n+j}(z^2)\big(2(z+1)\big)^{m-n-1-2j}.
\end{eqnarray}
\end{theorem}

\pf Let us first prove the recurrence relation \eqref{eqn 3.6}. Note that
$$N_{n,k}(z)=Res_v (f(v))^{k+1} v^{-n-2},$$
where
$f(v)$ is given by \eqref{eq-ff}, and $Res_v h(v)$ is the residue of $h(v)$ at $v=0$. Then we have
\begin{eqnarray}
F_{m,n}(z) \hskip-.22cm &=&\hskip-.22cm \lefteqn{
\sum_{k=0}^{n}(-z)^k\det\left(\begin{array}{cc}
N_{n,k}(z)     & N_{n,k+1}(z)      \\[5pt]
N_{m,k}(z)     & N_{m,k+1}(z)
\end{array}\right)}         \nonumber\\
\hskip-.22cm &=& \hskip-.22cm  \sum_{k=0}^{n} (-z)^{k}\left(Res_v(f(v))^{k+1} v^{-n-2} Res_u (f(u))^{k+2} u^{-m-2}\right.  \nonumber \\
\hskip-.22cm & & \hskip-.22cm  -\left.Res_v(f(v))^{k+2} v^{-n-2} Res_u (f(u))^{k+1} u^{-m-2} \right)                \nonumber   \\
\hskip-.22cm &=& \hskip-.22cm  Res_vRes_u \left(f(v)f(u)(f(u)-f(v)) v^{-n-2}u^{-m-2}\sum_{k=0}^{n} (-z)^{k}\left(f(v)f(u)\right)^{k}\right) \nonumber \\
\hskip-.22cm &=& \hskip-.22cm  Res_vRes_u \frac{f(v)f(u)(f(u)-f(v))}{1+zf(v)f(u)} v^{-n-2}u^{-m-2}. \label{eqn 3.8}
\end{eqnarray}

By (\ref{eqn 2.5}), it is clear that $f(v)$ has the compositional inverse $h(v)=\frac{v}{(1+v)(1+zv)}$.
Replacing $u$ by $h(u)$ and $v$ by $h(v)$ in (\ref{eqn 3.8}), we have
\begin{eqnarray*}
F_{m,n}(z) \hskip-.22cm &=&\hskip-.22cm  Res_vRes_u
\frac{vu(v-u)}{1+zuv} \frac{1}{h(v)^{n+2}h(u)^{m+2}}
                                     \frac{\partial h(v)}{\partial v}  \frac{\partial h(u)}{\partial u} \\
     \hskip-.22cm &=& \hskip-.22cm   Res_vRes_u \frac{(v-u)(1-zv^2)(1-zu^2)(1+v)^n(1+zv)^n(1+u)^m(1+zu)^m}{(1+zuv) v^{n+1}u^{m+1}}.
\end{eqnarray*}
Let
$$f_{m, n}(u,v,z)=\frac{(v-u)(1-zv^2)(1-zu^2)(1+v)^n(1+zv)^n(1+u)^m(1+zu)^m}{(1+zuv) v^{n+1}u^{m+1}}.$$
Thus, we have
\begin{align}\label{eq-key}
F_{m,n}(z)= Res_vRes_u f_{m, n}(u,v,z).
\end{align}
It is trivial to verify $F_{n,n}(z)=0$ (or $F_{m,n}(z)=-F_{n,m}(z)$)
by \eqref{eq-key}.

In order to prove \eqref{eqn 3.6}, let
$$ h_{m,n}(u,v,z) = f_{m,n}(u,v,z) - 2 (1+z) f_{m-1,n}(u,v,z) + f_{m-1,n+1}(u,v,z).$$
It suffices to check that $Res_vRes_u h_{m,n}(u,v,z)=0$.
By direct computation, we see that
\begin{multline*}
h_{m,n}(u,v,z)= \left(1+\frac{1}{u} \right)^{m-1}  (1+z u)^{m-1} \frac{(v^2-u^2)(1- z u^2)}{u^2}
\left(1+\frac{1}{v} \right)^n (1+z v)^{n} \frac{(1-z v^2)}{v^2}.
\end{multline*}
Now taking the residue respect to $u$, we obtain
\begin{multline*}
Res_u h_{m,n}(u,v,z) \\
= \left( 1+\frac{1}{v} \right)^n (1+z v)^{n} \frac{(1-z v^2)}{v^2}
\cdot \left(v^2 Res_u \left( 1+\frac{1}{u} \right)^{m-1}  (1+z u)^{m-1} \frac{(1- z u^2)}{u^2}  \right. \\ \left. - z Res_u \left(1+\frac{1}{u} \right)^{m-1}  (1+z u)^{m-1} (1- z u^2) \right).
\end{multline*}
Noting that
\begin{multline*}
Res_u \left( 1+\frac{1}{u} \right)^{m-1}  (1+z u)^{m-1} \frac{(1- z u^2)}{u^2} \\
= \sum_{k=0}^{m-1} {m-1 \choose k}{m-1 \choose k+1} z^{k+1} - z
   \sum_{k=1}^{m-1} {m-1 \choose k}{m-1 \choose k-1} z^{k-1}
=0,
\end{multline*}
we have
\begin{multline*}
Res_u h_{m,n}(u,v,z) = - z \left( 1+\frac{1}{v} \right)^n (1+z v)^{n} \frac{(1-z v^2)}{v^2} \\
\cdot Res_u \left(1+\frac{1}{u} \right)^{m-1}  (1+z u)^{m-1} (1- z u^2).
\end{multline*}
Now taking the residue with respect to $v$, we obtain
\begin{multline*}
Res_v Res_u h_{m,n}(u,v,z) = - z Res_v \left( 1+\frac{1}{v} \right)^n (1+z v)^{n} \frac{(1-z v^2)}{v^2} \\
 \cdot Res_u \left(1+\frac{1}{u} \right)^{m-1}  (1+z u)^{m-1} (1- z u^2) =0.
\end{multline*}

Hence, the recurrence (\ref{eqn 3.6}) is true.

We proceed to show that $F_{n+1,n}(z)=N_n(z^2)$ with the help of the
mathematica package {\tt HolonomicFunctions} \cite{H-Fun}. This
could be done along the following lines.

\begin{mma}
\In <<|RISC| ~\grave{} |HolonomicFunctions|~\grave{};\\
\end{mma}
\begin{mma}
\In  fmn :=\frac{(v-u)(1-zv^2)(1-zu^2)(1+v)^n(1+zv)^n(1+u)^m(1+zu)^m}{(1+zuv) v^{n+1}u^{m+1}};\\
\end{mma}
\begin{mma}
\In  fn :=fmn/. m -> n + 1;\\
\end{mma}
\begin{mma}
\In  |res| :=|FindCreativeTelescoping|[fn, \{Der[u], Der[v]\}, {S[n]}];\\
\end{mma}
\begin{mma}
\In  |L| :=|res|[[1]][[1]]\\
\end{mma}
\begin{mma}
\Out  (5 + n) (9 + 2 n) \mathbf{S}_n^3 + \big( 4 (2 + n) (4 + n) z - (z^2+1)(79 + 44 n + 6 n^2) \big)  \mathbf{S}_n^2
+ (2 + n) \big( (19 + 6 n) (z^4+1) -4 (5 + 2 n) (z^3+z) + 2 (1 + 2 n) z^2 \big) \mathbf{S}_n -2 (1 + n) (2 + n) (-1 + z)^4 (1 + z)^2\\
\end{mma}
By simplification we find that
\begin{align}\label{eq-rec-holo}
\mathbf{L }f_{n+1,n}(u,v,z) = \frac{\partial p_n(u,v,z)}{\partial u} + \frac{\partial q_n(u,v,z)}{\partial v},
\end{align}
where $p_n,q_n$ are two rational functions and $\mathbf{L}$ is the linear operator
\begin{multline*}
\mathbf{L }= (5 + n) (9 + 2 n) \mathbf{S}_n^3 + \big( 4 (2 + n) (4 + n) z - (z^2+1)(79 + 44 n + 6 n^2) \big)  \mathbf{S}_n^2 \\[5pt]
+ (2 + n) \big( (19 + 6 n) (z^4+1) -4 (5 + 2 n) (z^3+z) + 2 (1 + 2 n) z^2 \big) \mathbf{S}_n \\[5pt]
-2 (1 + n) (2 + n) (-1 + z)^4 (1 + z)^2
\end{multline*}
with $\mathbf{S}_n$ being the shift operator with respect to $n$. The explicit expressions of $p_n$ and $q_n$ are omitted. By \eqref{eq-rec-holo}, we get that
\[
\mathbf{L }\left(Res_v Res_u f_{n+1,n}(u,v,z)\right) = 0.
\]
On the other hand, it is easy to find that
\[
\mathbf{L}_1 N_n(z^2) = 0,
\]
where
\[
\mathbf{L}_1 = (4 + n) \mathbf{S}_n^2 - (5 + 2 n)(1+z^2) \mathbf{S}_n + (1 + n) (-1 + z)^2 (1 + z)^2.
\]
This could be done along the following lines.
\begin{mma}
\In |Nn| := \frac{1}{n+1}\binom{n+1}{i}\binom{n+1}{i+1}z^i;\\
\end{mma}
\begin{mma}
\In |res| := |CreativeTelescoping|[|Nn| /. z -> z^2, {S[i] - 1}, {S[n]}];\\
\end{mma}
\begin{mma}
\In |L1| := |res|[[1]][[1]]\\
\end{mma}
\begin{mma}
\Out (4 + n) \mathbf{S}_n^2 - (5 + 2 n)(1+z^2) \mathbf{S}_n + (1 + n) (-1 + z)^2 (1 + z)^2\\
\end{mma}
We have
\[
\mathbf{L} = \big( (9+2n) \mathbf{S}_n - 2 (2 + n) (-1 + z)^2 \big) \mathbf{L}_1,
\]
which could be easily implemented by using the following commands.
\begin{mma}
\In |ToOrePolynomial|[|L|, |OreAlgebra|[S[n]]];\\
\end{mma}
\begin{mma}
\In |ToOrePolynomial|[|L1|, |OreAlgebra|[S[n]]];\\
\end{mma}
\begin{mma}
\In |OreReduce|[|L|, \{|L1|\}, Extended -> True]\\
\end{mma}
\begin{mma}
\Out \{0, 1, \{\big( (9+2n) \mathbf{S}_n - 2 (2 + n) (-1 + z)^2 \big)\}\}\\
\end{mma}

Therefore, we obtain that
\[
Res_v Res_u f_{n+1,n}(u,v,z) = N_n(z^2)
\]
by checking the initial values of $n=0,1,2$. Thus, we have
$F_{n+1,n}(z)=N_n(z^2)$.

Finally, it is routine to verify that the polynomial sequence in the right of (\ref{eqn 3.7}) also satisfies the recurrence relation (\ref{eqn 3.6}) with the same initial conditions.
This completes the proof. \qed

When $z=1$ in Theorem \ref{theo 3.6}, routine simplification leads to the following result, which generalizes Theorem 4.4 of \cite{SunLu}.

\begin{corollary}
For any integers $m> n\geq 0$, there holds
\begin{eqnarray}\label{eqn 3.9}
\lefteqn{\sum_{k=0}^{n}(-1)^{k}\frac{(m-n)(k+1)(k+2)(2k+3)}{(n+1)(2n+3)(m+1)(2m+3)}\binom{2n+3}{n-k}\binom{2m+3}{m-k}  }\\
&=& \sum_{j=0}^{[\frac{m-n-1}{2}]}(-1)^j\binom{m-n-1-j}{j}C_{n+j+1} 4^{m-n-1-2j},  \nonumber
\end{eqnarray}
where $C_{i+1}=N_i(1)=\frac{1}{2i+3}\binom{2i+3}{i+1}$ is the $(i+1)$-th Catalan number \cite[A000108]{Sloane}.
\end{corollary}

\begin{remark}
Taking $n=0$ in (\ref{eqn 3.5}) and (\ref{eqn 3.7}), we get that  \begin{eqnarray*}
N_{m,1}(z)=\sum_{j=0}^{[\frac{m-1}{2}]}(-1)^j\binom{m-1-j}{j}N_{j}(z^2)\big(2(z+1)\big)^{m-1-2j}, \quad \mbox{for $m\geq 1$}.
\end{eqnarray*}
Using the recurrence for $N_{m,k}(z)$, i.e.,
\begin{align*}
N_{m,0}(z)&=(z+1)N_{m-1,0}(z)+zN_{m-1,1}(z)\ (n\geq 1),\\
N_{m,k}(z)&=N_{m-1,k-1}(z)+(z+1)N_{m-1,k}(z)+zN_{m-1,k+1}(z), \ (m,k\geq 1),
\end{align*}
by induction on $m$ and $k$ we have that each $N_{m,k}(z)$ can be represented as a rational combination of $N_{j}(z^2)(z+1)^{m-k-2j}$ for $0\leq j\leq [\frac{m-k}{2}]$, i.e.,
\begin{eqnarray*}
N_{m,k}(z)=\sum_{j=0}^{[\frac{m-k}{2}]}a_{m,k,j}N_{j}(z^2)(z+1)^{m-k-2j}
\end{eqnarray*}
for $a_{m,k,j}\in \mathbb{Q}$.
Thus it's natural to ask whether there exists any simple explicit expression or recurrence for $a_{m,k,j}$?
\end{remark}

\begin{remark}
According to (\ref{eqn 3.7}), $F_{m,n}(z)$ can be represented as a polynomial in $z$ of degree $m+n-1$ with symmetry coefficients for $m>n\geq 0$, i.e.,
\begin{eqnarray*}
F_{m,n}(z)=\sum_{j=0}^{m+n-1}b_{m,n,j}z^{j}, \ with\ b_{m,n,j}=b_{m,n,m+n-j-1}.
\end{eqnarray*}
Moreover, $F_{m,n}(z)$ can be represented as a polynomial with gamma basis $\{z^{j}(1+z)^{m+n-1-2j}\}_{j\geq 0}$ \cite{HSun}, that is,
\begin{eqnarray*}
F_{m,n}(z)=\sum_{j=0}^{[\frac{m+n-1}{2}]}c_{m,n,j}z^{j}(1+z)^{m+n-1-2j}.
\end{eqnarray*}
Does there exist any simple explicit expressions, recurrences or combinatorial interpretations for $b_{m,n,j}$ and $c_{m,n,j}$?
\end{remark}

\section{The recursive matrix $\mathcal{A}^{(x,y,y,\ldots),(0,0,0,\ldots)}$}\label{sect-4}

The aim of this section is to study certain weighted sums of $2\times 2$ minors of the recursive matrix $\mathcal{A}^{(x,y,y,\ldots),(0,0,0,\ldots)}$. Here we use
$M_{n,k}$ to denote the entry of $\mathcal{A}^{(x,y,y,\ldots),(0,0,0,\ldots)}$ in the $n$-th row and $k$-th column.

As in Section \ref{sect-3}, we first give an explicit formula for $M_{n,k}$.

\begin{lemma} For any $n\geq k\geq 0$, we have
\begin{align}\label{eq-mnk}
M_{n,k}&=\sum_{j=0}^{n-k}\binom{j+k-1}{j}x^{n-k-j}y^{j}.
\end{align}
\end{lemma}

\pf
By Lemma \ref{lem-connection}, the recursive matrix $\mathcal{A}^{(x,y,y,\ldots),(0,0,0,\ldots)}$ is the Riordan array $(g(v), f(v))$ with
\begin{align}
f(v)&=\frac{v}{1-yv}, \label{eq-f-4}\\
g(v)&=\frac{1}{1-xv}.\label{eq-g-4}
\end{align}
Thus
\begin{align*}
M_{n,k}&=[v^n]\frac{1}{1-xv}\cdot \frac{v^k}{(1-yv)^k}\\
&=\sum_{j=0}^{n-k}\binom{j+k-1}{j}x^{n-k-j}y^{j}.
\end{align*}
This completes the proof.
\qed

The few small values of $M_{n,k}$ are illustrated in Table 3.1.
\begin{center}
\begin{eqnarray*}
\begin{array}{|c|l|l|l|l|l|l}\hline
n/k & 0 & 1 & 2 & 3 & 4 \\\hline
0 & 1          &                          &                     &                &        \\
1 & x          & 1                        &                     &                &         \\
2 & x^2        & x+y                      & 1                   &                &         \\
3 & x^3        & x^2+xy+y^2               & x+2y                & 1              &         \\
4 & x^4        & x^3+x^2y+xy^2+y^3        & x^2+2xy+3y^2        & x+3y           & 1        \\\hline
\end{array}
\end{eqnarray*}
Table 3.1. The first values of $M_{n,k}$.
\end{center}

The main result of this section is as follows.

\begin{theorem} \label{thm-sect4} For any integers $m, n\geq 0$, there holds
\begin{eqnarray}
\lefteqn{ \sum_{k=0}^{\min\{m,n\}}y^{2k}\det\left(
\begin{array}{cc}
M_{n,k}     & M_{n,k+1}      \\[5pt]
M_{m,k}     & M_{m,k+1}
\end{array}\right) }   \nonumber\\
\hskip-.22cm &=&\hskip-.22cm  \sum_{k=1}^{\max\{m,n\}}\left(\binom{m+n-k}{n}-\binom{m+n-k}{m}\right)x^{k-1}y^{m+n-k}, \end{eqnarray}
where $M_{n,k}$ is given by \eqref{eq-mnk}.
\end{theorem}

\pf Note that $M_{n,k}=Res_v h(v)v^{-n-1}$, where
$$h(v)=\frac{1}{1-xv}\cdot \frac{v^k}{(1-yv)^k}.$$
Then we have
\begin{eqnarray*}
\lefteqn{\sum_{k=0}^{\min\{m,n\}}y^{2k}\det\left(
\begin{array}{cc}
M_{n,k}     & M_{n,k+1}      \\[5pt]
M_{m,k}    & M_{m,k+1}
\end{array}\right)}  \\
\hskip-.22cm &=& \hskip-.22cm  \sum_{k=0}^{\min\{m,n\}}y^{2k}\left(Res_v h(v)v^{-n-1} Res_u h(u)u^{-m-1}\right.  \\
\hskip-.22cm & & \hskip-.22cm  -\left.Res_v h(v)v^{-n-1} Res_u h(u)u^{-m-1} \right)                             \\
\hskip-.22cm &=& \hskip-.22cm  \sum_{k=0}^{\min\{m,n\}}y^{2k}\left(Res_v \frac{1}{1-xv}\left(\frac{v}{1-yv}\right)^{k}v^{-n-1}Res_u\frac{1}{1-xu}\left(\frac{u}{1-yu}\right)^{k+1}u^{-m-1} \right.   \\
\hskip-.22cm & & \hskip-.22cm  -\left. Res_v \frac{1}{1-xv}\left(\frac{v}{1-yv}\right)^{k+1}v^{-n-1}Res_u\frac{1}{1-xu}\left(\frac{u}{1-yu}\right)^{k}u^{-m-1} \right)                             \\
\hskip-.22cm &=& \hskip-.22cm  Res_vRes_u\frac{1}{1-xv}\frac{1}{1-xu}v^{-n-1}u^{-m-1}\left(\frac{u}{1-yu}-\frac{v}{1-yv}\right) \sum_{k=0}^{\min\{m,n\}} \left(\frac{yv}{1-yv}\right)^{k} \left(\frac{yu}{1-yu}\right)^{k}   \\
\hskip-.22cm &=& \hskip-.22cm  Res_vRes_u\frac{1}{1-xv}\frac{1}{1-xu}v^{-n-1}u^{-m-1}\frac{u-v}{(1-yu)(1-yv)} \frac{1}{1-\frac{yv}{1-yv}\frac{yu}{1-yu}}   \\
\hskip-.22cm &=& \hskip-.22cm  Res_vRes_u x^{-1}\left(\frac{1}{1-xu}-\frac{1}{1-xv}\right)\frac{1}{1-y(u+v)} v^{-n-1}u^{-m-1}  \\
\hskip-.22cm &=& \hskip-.22cm  Res_vRes_u x^{-1}\sum_{k=0}^{\infty}x^{k}u^{k}\sum_{i=0}^{\infty}y^{i}\sum_{j=0}^{i}\binom{i}{j}v^{j}u^{i-j}v^{-n-1}u^{-m-1}   \\
\hskip-.22cm & & \hskip-.22cm  -Res_vRes_u x^{-1}\sum_{k=0}^{\infty}x^{k}v^{k}\sum_{i=0}^{\infty}y^{i}\sum_{j=0}^{i}\binom{i}{j}u^{j}v^{i-j}v^{-n-1}u^{-m-1} \\
\hskip-.22cm &=& \hskip-.22cm  \sum_{k=0}^{n}\binom{m+n-k}{n}x^{k-1}y^{m+n-k}-\sum_{k=0}^{m}\binom{m+n-k}{m}x^{k-1}y^{m+n-k}   \\
\hskip-.22cm &=& \hskip-.22cm  \sum_{k=1}^{\max\{m,n\}}\left(\binom{m+n-k}{n}-\binom{m+n-k}{m}\right)x^{k-1}y^{m+n-k} .
\end{eqnarray*}
This completes the proof. \qed

Letting $m=n+1$ in Theorem \ref{thm-sect4}, we obtain the following result related to the ballot numbers.

\begin{corollary} For any integer $n\geq 0$, there holds
\begin{eqnarray*}
\sum_{k=0}^{n}y^{2k}\det\left(
\begin{array}{cc}
M_{n,k}       & M_{n,k+1}     \\[5pt]
M_{n+1,k}     & M_{n+1,k+1}
\end{array}\right)
\hskip-.22cm &=&\hskip-.22cm  \sum_{k=0}^{n}C_{n,k}x^ky^{2n-k}.
\end{eqnarray*}
where $M_{n,k}$ is given by \eqref{eq-mnk} and $C_{n,k}=\frac{k+1}{n+1}\binom{2n-k}{n}$ are the ballot numbers \cite[A033184]{Sloane}.
\end{corollary}

\section*{Acknowledgements}

The first author and the third author are supported in part by ``Liaoning BaiQianWan Talents Program".
The second author is supported in part by the National Science Foundation of China (No. 11771330).
The fourth author is supported in part by the National Science Foundation of China (Nos. 11231004, 11522110).

%==============================================================================================================

\end{document}